И. В. Баяк

# О порождении линейных групп комбинаторными группами

*Резюме*: Описана группа стрелочных подстановок и процедура генерации полной линейной группы и некоторых ее подгрупп.

1. **Группа стрелочных подстановок как расширение симметрической группы**

   В комбинаторных построениях данного пункта без объяснения используется понятие цикла [1] и системы порождающих [2]. В остальном, изложение логически замкнуто и не требует дополнительных ссылок. Итак, пусть $I = \{1,..,n\}$, $s : I \to I$ - биекция а $s's : I \to I (s's(i) = s'(s(i)))$ - композиция биекций. Тогда множество всех $n!$ биекций $S = \{s\}$ относительно произведения $S \times S \to S : (s', s) \to s's$ составляет группу подстановок (симметрическую группу) степени $n$. Вместе с тем, множество всех четных биекций составляет подгруппу $S^+$ группы $S$, которая в общем случае порождается 3-циклами, а именно, $S^+(n) = \langle \{(i, j, k)\}_I \rangle$, где $n \geq 3$ и $i \neq j \neq k$, но если $n = 1,2$, то $S^+(1) = S^+(2) = e$, где $e$ - тождественная подстановка.

   Пусть также $A = \{\pm 1,..,\pm n\}$ а $p : I \to A$ - такая инъекция, что $|p| : I \to I$ - биекция, причем существует композиция этих инъекций $p'p : I \to A$: $p'p(i) = \mathrm{sgn}(p(i)) p'(|p(i)|)$. Тогда множество всех $2^n n!$ биективных по модулю инъекций $P = \{p\}$ относительно определенной ранее композиции, т.е. произведения $P \times P \to P : (p', p) \to p'p$, составляет группу, называемую нами группой стрелочных подстановок n-й степени. Вместе с тем, стрелочная подстановка называется стрелочной транспозицией, если это элементарная перестановка, т.е. $p^{j+1}(l) = p^j(m)$, $p^{j+1}(m) = p^j(l)$, $p^{j+1}(i) = p^j(i)$ для $l, m, i \in I \wedge i \neq l, m$, или элементарная инверсия, т.е. $p^{j+1}(k) = -p^j(k)$, $p^{j+1}(i) = p^j(i)$ для $k, i \in I \wedge i \neq k$, и где $j$ - рекурсивный индекс, $i$ - бегущий индекс.

   По аналогии с простыми подстановками, всякая стрелочная подстановка может быть получена композицией стрелочных транспозиций, а четность стрелочной подстановки $p$ определяется четностью числа стрелочных транспозиций для перехода к $p$. При этом, четность подстановки $p$ инвариантна относительно выбора композиции транспозиций для перехода к $p$, т.е. $(-1)^{\sigma(p)} = (-1)^{\sigma_1(p)} = (-1)^{\sigma_2(p)}$, где $\sigma_1(p)$ и $\sigma_2(p)$ - количество транспозиций в композиции 1 и 2 соответственно, а $\sigma(p) \in Z/2Z$. Итак, если $(-1)^{\sigma(p)} = 1$, то $p$ - четная подстановка, если же $(-1)^{\sigma(p)} = -1$, то $p$ - нечетная подстановка. Множество всех четных стрелочных подстановок составляет подгруппу $P^+$ группы $P$, которая в общем случае порождается генераторами

2$(j,-k)$: $j \to (-k), k \to j, i \to i$, где $i \neq j \neq k$, а именно, $P^+(n) = \langle\{(j,-k)\}\rangle$, но если $n = 1$, то $P^+(1) = e$ - тождественная стрелочная подстановка.

Далее пусть $I_1 = \{1,..,k\}$, $I_2 = \{k+1,..,n\}$, $I_1 \cup I_2 = I$, $s_1 = s(I_1)$, $s_2 = s(I_2)$. Тогда $s$ можно представить как составную подстановку $s_1 \times s_2$, состоящую из подстановки по месту $I_1$ и подстановки по месту $I_2$, т.е. из размещений $s_1$ и $s_2$. Вместе с тем, размещение $s_1(s_2)$ называется транспозицией по месту $I_1(I_2)$, если это элементарная перестановка внутри $I_1(I_2)$, или элементарное замещение, а именно, для замещения из $I_2$ в $I_1$ имеем $s_1^{j+1}(l) = s^j(m)$, $s_1^{j+1}(i) = s_1^j(i)$ где $m \in I_2$; $l, i \in I_1 \wedge i \neq l$, а для замещения из $I_1$ в $I_2$ имеем $s_2^{j+1}(l) = s^j(m)$, $s_2^{j+1}(i) = s_2^j(i)$ где $m \in I_1$; $l, i \in I_2 \wedge i \neq l$. Четность составной подстановки $s_1 \times s_2$ задается четностью ее размещений, т.е. числом $(-1)^{\sigma(s_1)} \cdot (-1)^{\sigma(s_2)}$, где $\sigma(s_1)$, $\sigma(s_2)$ - факторизованное число транспозиций по месту $I_1$, $I_2$ для перехода к размещению $s_1$, $s_2$ соответственно. При этом, если $(-1)^{\sigma(s_1)} \cdot (-1)^{\sigma(s_2)} = 1$, то $s_1 \times s_2$ называется четно-составной подстановкой, а множество всех четно-составных подстановок составляет подгруппу $S^+(k, n-k)$ группы $S(n)$, которая в общем случае порождается 3-циклами, действующими внутри подмножеств, и произвольным 2-циклом, действующим между $I_1$ и $I_2$, а именно, $S^+(k, n-k) = \langle \{(i,j,o)\}_{I_1}, (l,m), \{(i,j,o)\}_{I_2} \rangle$, где $l \in I_1$ а $m \in I_2$.

В свою очередь, поскольку понятие составной подстановки допускает естественное расширение в группу стрелочных подстановок, то множество всех четно-составных стрелочных подстановок составляет подгруппу $P^+(k, n-k)$ группы $P(n)$, которая в общем случае порождается генераторами $(i,-j)$, действующими внутри подмножеств, и произвольной парой генераторов $\pm(l,m)$, действующей между ними, где $+(l,m)$ - это элементарная перестановка, а $-(l,m)$ - это перестановка с инверсиями, т.е. $-(l,m): l \to (-m), m \to (-l)$. Тем самым, имеем $P^+(k, n-k) = \langle \{(i,-j)\}_{I_1}, \pm(l,m), \{(i,-j)\}_{I_2} \rangle$, где $l \in I_1$ а $m \in I_2$.

Далее, прежде чем приступить к матричной реализации группы стрелочных подстановок, определим детерминант подмножества строк квадратной матрицы. Пусть $I^* \subset I$; $b: I \to I$ - биекция, тогда $b(I^*): I^* \to I$ - инъекция, т.е. подстановка по месту $I^*$ или размещение. Если $b^*: I^* \to I$ - произвольная инъекция, то $\{b(I^*)\} \approx \{b^*\}$ и четность размещения $b^*$ определяется числом $(-1)^{\sigma(b^*)}$, где $\sigma(b^*)$ - количество транспозиций для перехода от $e^*$ к $b^*$, где $e^*$ - тождественное отображение. Тем самым, если $A = (a_{ij})_{i,j \in I}$, $A(I^*) = (a_{ij})_{i \in I^*, j \in I}$, тогда принимаем, что $\det(A(I^*)) = \sum_{\{b^*\}} \prod_{i \in I^*} a_{ib^*(i)} (-1)^{\sigma(b^*)}$.

Пусть теперь $\{(a_{ij})_I\}$ - множество всех переходных матриц, т.е. таких квадратных матриц, в которых каждый столбец и каждая строка имеют один и только один ненулевой элемент, причем, если $a_{ij} \neq 0$, то $a_{ij} \in \{\pm 1\}$. Тогда для группы стрелочных подстановок $P(n)$ имеем изоморфизм $P(n) \approx \left\{ \left( a_{ij} \right)_{i,j \in I} \right\}$, где $a_{ij} = \text{sgn}(p(i)) \ \forall \ j = |p(i)|$ и $a_{ij} = 0 \ \forall \ j \neq |p(i)|$. Кроме того, $P^+(I) \approx \left\{ \left( a_{ij} \right)_I \middle| \det A = +1 \right\}$ а $P^+(I^*, I \setminus I^*) \approx \left\{ \left( a_{ij} \right)_I \middle| \det A(I^*) \det A(I \setminus I^*) = +1 \right\}$, где $P^+(I)$, $P^+(I^*, I \setminus I^*)$ - подгруппы четных и четно-составных стрелочных подстановок соответственно. Вместе с тем, если $J = \{1,..,m\}$; $m \leq n$; $I_J \equiv \left\{ I_j \right\}_J$, где $\bigcup_m I_j = I, \bigcap_m I_j = \emptyset, \text{card}(I_j) = n_j, \sum_m n_j = n$, причем подсемейство $I_{j+1}$ заполняется последовательной выборкой из $I$ вслед за заполнением $I_j$, а $I_1 = \{1,..,n_1\}$, то $P^+(I_J) \equiv P^+(n_1,...,n_m) \approx \left\{ \left( a_{ij} \right) \middle| \det A(I_1) \cdots \det A(I_m) = +1 \right\}$.

## 2. Группа стрелочных подстановок как генератор полной линейной группы

Пусть $R_+ \equiv \{x \in R | x > 0\}$; $RP(n)$ - групповая алгебра над $R$ с базисом $P(n)$, которая порождается как линейная оболочка, натянутая на группу $P(n)$, т.е. $RP(n) = \langle P(n) \rangle$; $GP(n)$ - мультипликативная группа алгебры $RP(n)$; а $\text{Aut } RP(n) \equiv GP(n)/R_+$ - группа автоморфизмов этой же алгебры. Тогда имеет место основное утверждение.

**Lemma 1:** $RP(n) \approx \text{End } R^n \equiv ML(n), \ GP(n) \approx GL_n R \equiv GL(n),$
$$\text{Aut } RP(n) \approx \langle SL(n), P(n) \rangle$$

Действительно, поскольку $P(n)$ реализуется множеством переходных матриц $A(I)$, из которых всегда можно выбрать $n^2$ линейно независимых матриц, то $A(I)$ включает в себя и некоторый базис пространства $ML(n)$, имеющего размерность $n^2$, а следовательно линейная оболочка, натянутая на $A(I)$, совпадает с $ML(n)$. В свою очередь, поскольку $GL(n) = \langle SL(n), A(I), R_+ \cdot E \rangle$, где $E$ - единичная матрица, то $GL(n)/R_+ \approx \langle SL(n), P(n) \rangle$ и лемма доказана.

Дополнительным источником процедуры порождения линейных групп служит утверждение о генерации общей линейной группы.

**Lemma 2:**

Пусть $GL_{i,j} \equiv \text{diag}\left[ 1,..,GL(2)_{i,j}^{i,j},..,1 \right]_n$, где $i \neq j$, причем верхние индексы указывают номера строк, а нижние - номера столбцов, на пересечении которых располагается группа $GL(2)$. Тогда
$GL(n) = \left\langle \left\{ GL_{i,j} \right\}_I \right\rangle = \left\langle \left\{ GL_{i,i+1} \right\}_I \right\rangle$.

Действительно, поскольку $P(n) = \left\langle \left\{ P(2)_{i,j} \right\}_I \right\rangle = \left\langle \left\{ P(2)_{i,i+1} \right\}_I \right\rangle$, где $P(2)_{i,j}$ - подгруппа стрелочных подстановок степени $n$, изоморфная группе $P(2)$, которая представлена множеством стрелочных подстановок по месту $i, j$, то $GP(n) = \langle \{GP(2)_{i,j}\}_I \rangle$, но $GP(2)_{i,j} \approx GL_{i,j}$, а следовательно лемма доказана.

Далее пусть даны циклические группы стрелочных подстановок второй степени $P^+(2) \approx \left\langle \begin{pmatrix} 0 & 1 \\ -1 & 0 \end{pmatrix} \right\rangle$, $P^+(1,1) \approx \left\langle \begin{pmatrix} 0 & 1 \\ 1 & 0 \end{pmatrix} \right\rangle$, $P^-(2) \approx \left\langle \begin{pmatrix} 1 & 0 \\ 0 & -1 \end{pmatrix} \right\rangle$, которые вместе с тривиальной подгруппой исчерпывают все подгруппы группы $P(2)$. Тогда после элементарных построений имеем изоморфизмы $Aut\, RP^+(2) \approx \left\{ \begin{pmatrix} \cos x & \sin x \\ -\sin x & \cos x \end{pmatrix}, x \in R \right\} \equiv SO(2)$, $Aut\, RP^+(1,1) \approx \left\{ \pm \begin{pmatrix} chx & shx \\ shx & chx \end{pmatrix}, x \in R \right\} \equiv SO(1,1)$, $Aut\, RP^-(2) \approx \left\{ \pm \begin{pmatrix} e^x & 0 \\ 0 & e^{-x} \end{pmatrix}, x \in R \right\} \equiv ST(2)$, а поскольку $P(2) = P^+(2) \oplus P^+(1,1) \oplus P^-(2)$, то $GL(2) = \langle SO(2), SO(1,1), ST(2), P(2), R_+ \rangle$ и поэтому $SL(2) = \langle SO(2), SO(1,1), ST(2) \rangle$. Пусть теперь $SO(2)_{i,j} \equiv diag[1,...,SO(2)_{i,j}^{i,j},...,1]_n$, тогда $Aut\, RP^+(n) \approx \langle \{SO(2)_{i,i+1}\}_I \rangle = SO(n)$. Аналогично, если $SO(1,1)_{i,j} \equiv diag[1,...,SO(1,1)_{i,j}^{i,j},...,1]_n$, то $Aut\, RP^+(m, n-m) \approx \langle \{SO(2)_{i,i+1}\}_{I_1}, SO(1,1)_{m,m+1}, \{SO(2)_{i,i+1}\}_{I_2} \rangle = SO(m, n-m)$. Более того, возможно расширение специальной ортогональной группы, а именно, $Aut\, RP^+(I_J) \approx \langle \{\{SO(2)_{i,i+1}\}_{I_j}\}_J, \{SO(1,1)_{n_j,n_j+1}\}_J \rangle \equiv SO(I_J)$. Если же $SL_{i,j} \equiv diag[1,...,SL(2)_{i,j}^{i,j},...,1]_n$, то $Aut^+RP(n) = \langle \{SL_{i,i+1}\}_I \rangle = SL(n)$, где $Aut^+RP(n)$ - автоморфизмы алгебры, сохраняющие ее ориентацию, т.е. $Aut\, RP(n) = \langle Aut^+RP(n), P(n) \rangle$.

В целом понятно, что всевозможные подгруппы группы стрелочных подстановок порождают всевозможные линейные группы (в том числе – унитарные, симплектические и т.д.), и поэтому задача классификации, описания и конструирования линейных групп сводится к решению такой же задачи для подгрупп конечной группы.